\documentclass[12pt]{article}
\usepackage[reqno]{amsmath}
\usepackage{amssymb, amsthm, enumerate}
\usepackage{graphicx}
\usepackage{array}
\usepackage[svgnames]{xcolor}
\usepackage{tikz}
\usepackage{mathtools}
\usepackage{url}
\usetikzlibrary{decorations.markings}
\usetikzlibrary{shapes.geometric}

\usepackage{fullpage}

\pgfdeclarelayer{edgelayer}
\pgfdeclarelayer{nodelayer}
\pgfsetlayers{edgelayer,nodelayer,main}

\tikzstyle{none}=[inner sep=0pt]

\tikzstyle{rn}=[circle,fill=Red,draw=Black,line width=0.8 pt]
\tikzstyle{gn}=[circle,fill=Lime,draw=Black,line width=0.8 pt]
\tikzstyle{yn}=[circle,fill=Yellow,draw=Black,line width=0.8 pt]
\tikzstyle{blackcirc}=[circle,fill=Black,draw=Black, scale= 0.5]
\tikzstyle{whitecirc}=[circle,fill=White,draw=Black, scale=0.5]
\tikzstyle{newstyle}=[circle,fill=White,draw=Black]

\tikzstyle{simple}=[-,draw=Black,line width=1.000]
\tikzstyle{arrow}=[-,draw=Black,postaction={decorate},decoration={markings,mark=at position .5 with {\arrow{>}}},line width=2.000]
\tikzstyle{tick}=[-,draw=Black,postaction={decorate},decoration={markings,mark=at position .5 with {\draw (0,-0.1) -- (0,0.1);}},line width=2.000]

\expandafter\ifx\csname pdfpagebox\endcsname\relax\else
\fi

{ %\theorembodyfont{\slshape}
\newtheorem{theorem}{Theorem.}[section]
\newtheorem{lemma}[theorem]{Lemma.}
\newtheorem{corollary}[theorem]{Corollary.}
\newtheorem{proposition}[theorem]{Proposition.}}
\newtheorem*{thm:offdiag}{Theorem \ref{thm:offdiag}}

\def\proof{\noindent{{\sl Proof. }}}
\def\sqr#1#2{{\vbox{\hrule height.#2pt
    \hbox{\vrule width.#2pt height#1pt \kern#1pt
        \vrule width.#2pt}\hrule height.#2pt}}}
\def\eqed{\sqr53}
\def\qed{%
    \ifmmode\eqno\eqed
    \else\nobreak\ \hfill\eqed\medbreak\fi}

\newcommand\Zv{{\mathbf v}}

\newcommand\Zx{{\mathbf x}}

\newcommand\cx{{\mathbb C}}% complexes

\newcommand\re{{\mathbb R}}%reals

\newcommand\comp[1]{{\mkern2mu\overline{\mkern-2mu#1}}}

\newcommand\pmat[1]{\begin{pmatrix} #1 \end{pmatrix}}

\DeclareMathOperator\diag{diag}

\DeclareMathOperator\rea{Re}

\title{Spectral bound for separations in Eulerian digraphs}

\author{Krystal Guo\thanks{Supported in part by NSERC. Most of this work was done when the author was at Department of Mathematics, Simon Fraser University, Burnaby, B.C. V5A 1S6. } \\[1mm]
   {Department of Combinatorics \& Optimization}\\{University of Waterloo}\\{Waterloo, Ontario N2L 3G1} \\ \texttt{kguo@uwaterloo.ca}
 }

\begin{document}
\maketitle

\begin{abstract}

The spectra of digraphs, unlike those of graphs, is a relatively 
unexplored territory. In a digraph, a separation is a pair of sets of 
vertices $X$ and $Y$ such that there are no arcs from $X$ and $Y$. For 
a subclass of eulerian digraphs, we give an bound on the size of a 
separation in terms of the eigenvalues of the Laplacian matrix.

\vspace{5pt}
\noindent Keywords: algebraic graph theory, eigenvalue, directed graphs

\vspace{5pt}
\noindent Mathematical Subject Classification: 05C50, 05C20
\end{abstract}

\section{Introduction}

The theory of graph spectra is a rich and interesting area. There has been extensive study about the interplay of eigenvalues of a graph and various graph properties, such as the diameter \cite{C89, M91} or the chromatic number \cite{Cv72, H70, H95}; see also \cite{MoPo93}. The relationship between symmetries of a graph and its eigenvalues has also been investigated extensively, for example in \cite{CG97, R07, R13}. The eigenvalues of the Laplacian matrix of a graph determine the number of connected components. For connected graphs, the eigenvalues of the adjacency matrix determine whether the graph is bipartite. There are spectral bounds on the independence number and many other properties. 

For digraphs, in contrast, there are relatively few results. There is a directed analogue of Wilf's bound on chromatic number \cite{M10}, however other directed analogues are yet to be seen. The adjacency matrix of a digraph is usually difficult to work with. It is not always diagonalizable and may have complex eigenvalues. In addition, the interlacing theorem does not hold for adjacency matrices of digraphs, in general. The digraphs we consider here are also known as \textsl{mixed graphs} in the literature. 
%All acyclic digraphs have the same spectrum as the empty graph. 

The interlacing theorem is a powerful tool for studying with the eigenvalues of graphs. We would like to extend its usage to some classes of digraphs, with possibly different choices of matrices. 
To this end, we use the $A$-Laplacian matrix, $L(X) = \Delta^+(X) - A(X)$, where $\Delta^+(X)$ is the diagonal matrix with diagonal entries equal to the out-degrees of vertices of the digraph. We prove a result which generalizes the following theorem of Haemers for graphs to the class of digraphs with normal $A$-Laplacian matrices.

\begin{theorem}[\cite{H95}]\label{thm:haemers} Let $X$ be a connected graph on $n$ vertices and let $(Z,Y)$ be a separation in $X$. Then
\[
\frac{|Y| |Z|}{(n-|Y|)(n - |Z|)} \leq  \frac{ |\alpha + \sigma_n |^2}{\alpha^2},
\]
where $\alpha = -\frac{1}{2}(\sigma_2 + \sigma_n)$ and $0 = \sigma_1 < \sigma_2 \leq \cdots \leq \sigma_n$ are the Laplacian eigenvalues of $X$.
\end{theorem}

We give the following generalization of this theorem to digraphs whose $A$-Laplacian matrix is normal. 

\begin{thm:offdiag} 
 Let $X$ be a connected digraph on $n$ vertices where $L(X)$ is normal and $(Z, Y)$ be a separation in $X$. Then
\[
\frac{|Y| |Z|}{(n-|Y|)(n - |Z|)} \leq \frac{ |\alpha + \theta |^2}{\alpha^2} ,
\]
where $ \alpha = \begin{cases} -f(\theta) - f(\nu), & \text{if } \rea(\lambda) \geq \rea(\theta) \text{ for all }\lambda \notin \{0, \theta\}; \\
-f(\theta) - g(\mu), & \text{otherwise} \end{cases} $ 
and 
 \begin{itemize}
 \item  $ f(\lambda) = \frac{|\lambda|^2}{2 \rea(\lambda)} $,
 \item $\theta \neq 0$ is the eigenvalue of $L$ which maximizes $f$ amongst non-zero eigenvalues of $L$ 
 \item $\nu \neq 0$ is the eigenvalue of $L$ which minimizes $f$ amongst non-zero eigenvalues of $L$ 
 \item $ g(\lambda) = \frac{\rea(\lambda)(f(\theta) - f(\lambda))}{\rea(\theta) - \rea(\lambda)} $, and
 \item $\mu$ is the eigenvalue of $L(X)$ which minimizes $g$ such that $g(\mu)>0$, if such an eigenvalue exists, and $\mu =0$ otherwise.
 \end{itemize}
\end{thm:offdiag}

In Section \ref{sec:digrs}, we will give a combinatorial description of digraphs whose $A$-Laplacian matrices are normal. For example, any Cayley digraph of an abelian group will have a normal $A$-Laplacian matrix. In Section \ref{sec:eig}, we will prove some preliminary facts about the eigenvalues of digraphs with normal $A$-Laplacian. We will give the proof of the main theorem, Theorem \ref{thm:offdiag}, in Section \ref{sec:interlacing}, as well as a corollary for tournaments. 

\section{Directed graphs with normal Laplacian matrix}\label{sec:digrs}

We will define the $A$-Laplacian matrix of directed graphs as a directed analogy of the usual graph Laplacian. For a digraph $X$, let $\Delta^+(X)$ be the diagonal matrix indexed by the vertex set of $X$ with $\Delta^+(u,u)$ equal to the out-degree of vertex $u$. Let $A$ be the adjacency matrix of $X$. Then, the \textsl{$A$-Laplacian matrix of digraph} $X$ is 
\[
L(X) = \Delta^+(X) -A(X). 
\]
We will write $L$, $\Delta^+$ and $A$ when there is no ambiguity. We are interested in digraphs where $L$ is normal. We will assume our digraphs have no loops or parallel arcs. 

We can easily give a combinatorial description of digraphs whose $A$-Laplacian is normal. Let $X$ be a digraph. We denote by $d^+(u)$ the out-degree of a vertex $u$. Let $d^+(u,v)$ be the common out-neighbours of $u$ and $v$; that is 
\[
d^+(u,v) = |\{w \mid uw, vw \in E(X)\}|.
\]
Similarly, let $d^-(u,v) = |\{w \mid uw, vw \in E(X)\}|$.

\begin{lemma} Let $X$ be a digraph. The $A$-Laplacian of $X$ is normal if and only if, for every pair of vertices $u,v$,
\[
d^-(u,v) - d^+(u,v) = \begin{cases} 0 & \text{if } uv, vu \in E(X) \text{ or } uv, vu \notin E(X); \\
d(u)-d(v) & \text{if } uv\in E(X) \text{ and } vu \notin E(X); \text{ and,}\\
d(v)-d(u) & \text{if } vu\in E(X) \text{ and } uv \notin E(X).\\
\end{cases}
\]
\end{lemma}

\proof We will look at the $(u,v)$ entry of the matrices $L^TL$ and $LL^T$. We have that 
\[\begin{split}
(L^TL)(u,v) &= \sum_{w\in V(X)} L^T(u,w) L(w,v) \\
&= \sum_{w\in V(X)} L(w,u) L(w,v) \\
&= d^-(u,v) + d(u) L(u,v) + d(v)L(v,u).
\end{split}
\]
Similarly, 
\[\begin{split}
(LL^T)(u,v) &= \sum_{w\in V(X)} L(u,w) L^T(w,v) \\
&= \sum_{w\in V(X)} L(u,w) L(v,w) \\
&= d^+(u,v) + d(u) L(v,u) + d(v)L(u,v).
\end{split}
\]
Note that we always have $(L^TL)(u,u) = (LL^T)(u,u)$. For $u\neq v$, we have $(L^TL)(u,v) = (LL^T)(u,v)$ iff
\[
d^-(u,v) - d^+(u,v) = \begin{cases} 0 & \text{if } L(v,u) = L(u,v); \\
d(u)-d(v) & \text{if } L(u,v) = -1 \text{ and } L(v,u)=0; \text{ and,}\\
d(v)-d(u) & \text{if } L(v,u) = -1 \text{ and } L(u,v)=0.\\
\end{cases} \]
whence the lemma follows. \qed 

Furthermore, we can show the digraphs with normal Laplacian are a subclass of eulerian digraphs. 

\begin{proposition} \label{prop:eulerian} For a weakly connected digraph $X$, if $L(X)$ is normal, then $X$ is eulerian. \end{proposition}

\proof Since $L$ is normal, we have that 
\[ \begin{split}
0 &= LL^T - L^TL \\
&= (\Delta^+ - A)(\Delta^+ - A)^T - (\Delta^+ - A)^T(\Delta^+ - A) \\
&= (\Delta^+)^2 - A\Delta^+ - \Delta^+ A^T + AA^T - (\Delta^+)^2 + A^T\Delta^+ + \Delta^+ A - A^TA \\
&= \Delta^+(A - A^T) + (A^T - A)\Delta^+ + AA^T - A^TA =: M.
\end{split}
\]
For a square matrix $N$, let $\diag(N)$ denote the vector consisting of the diagonal entries of $N$:
\[
(\diag(N))_u = N(u,u).
\]
 Observe that $\Delta^+$ is a diagonal matrix and both $A$ and $A^T$ have zero diagonal. Then
\[
\diag(\Delta^+(A - A^T)) = \diag((A^T - A)\Delta^+) = 0.
\]
We see that $M=0$ and so has $\diag(M) = 0$. Then
\begin{equation}\label{eq:aat}
\diag(AA^T) = \diag(A^TA).
\end{equation}
Combinatorially, we have that
\[
(A^T A)_{u,u} = |\{w \,:\, wu \text{ is an arc of } X \}| = \deg^-(u)
\]
and
\[
(A A^T)_{u,u} = |\{y \,:\, uy \text{ is an arc of } X \}| = \deg^+(u).
\]
This and (\ref{eq:aat}) implies that $\deg^-(u) = \deg^+(u)$ for every vertex $u$ of $X$ and, since $X$ is weakly connected, so $X$ is eulerian. \qed

If a weakly connected digraph is Eulerian, then it is also strongly connected. Since there is no confusion, we may say that such graphs are $\textsl{connected}$.  

We observe that, if the digraph is regular, then the $A$-Laplacian matrix is normal if and only if the adjacency matrix is normal. We say that a digraph is \textsl{normal} if it has a normal Laplacian matrix and a normal adjacency matrix. 

\begin{theorem} Every Cayley digraph on an abelian group is normal. \end{theorem}

\proof Consider a  Cayley digraph $X = \text{Cay}(G, C)$ where $G$ is abelian and let $A$ be the adjacency matrix of $X$. Since $X$ is regular, we need only check that $A$ is normal. Then
\[
(A^T A)_{u,v} = |\{w \,:\, wu \text{ and } wv \text{ are arcs of } X \}|
\]
and
\[
(A A^T)_{u,v} = |\{y \,:\, uy \text{ and } vy \text{ are arcs of } X \}|.
\]
Let $u,v$ be vertices of $X$ and suppose they have a common in-neighbour $w$. Then for some $a_1 \neq a_2 \in C$, we have
\[
u = a_1w \text{ and } v = a_2 w. 
\]
It is clear that $a_2 u$ is an out-neighbour of $u$ and $a_1 v$ is an out-neighbour of $v$. In addition,
\[
a_2u = a_2 a_1 w = a_1 a_2 w = a_1v
\]
and so $a_2 u$ is a common out-neighbour of $u$ and $v$. Conversely, following the same argument, given a common out-neighbour of $u$ and $v$, we may construct a common in-neighbour of $u$ and $v$. Then
\[
(A^TA)_{u,v} = (AA^T)_{u,v}
\]
and so $A$ is a normal matrix. \qed

In general, the adjacency matrix being normal does not have to coincide with the Laplacian matrix being normal. We present some data on all digraph on 4,5 and 6 vertices. Note that if $X$ has a normal adjacency matrix or a normal Laplacian, it must be eulerian.

\begin{table}[htdp]
\caption{Small digraphs with normal Laplacian and adjacency matrices.}
\begin{center}
\begin{tabular}{|l|c|c|c|}
\hline
& 4 vertices & 5 vertices & 6 vertices \\
\hline
number of digraphs & 218 & 9608 & 1540944 \\
eulerian & 17 & 107 & 2269\\
regular & 5 & 10 & 52\\
normal Laplacian & 14 & 43 & 194\\
normal adjacency matrix & 14 & 45 & 212 \\
normal & 14 & 43 & 190 \\
connected and eulerian & 12 & 90 & 2162\\
undirected  & 10 & 31 & 43\\
\hline
\end{tabular}
\end{center}
\label{tab:sm}
\end{table}%

\section{Eigenvalues of normal Laplacian matrices}\label{sec:eig}

To prove the main results, we need the following technical lemma. Lemma \ref{lem:laplace} can be proved by considering $L + L^T$ as the Laplacian of a multi-graph and appealing to known results about Laplacians of graphs. Here, however, we will give a direct proof. Parts (a) and (b) of Lemma \ref{lem:laplace} follow from well-known results on $M$-matrices, see \cite{BePl94}.

\begin{lemma}\label{lem:laplace} If $X$ is a connected digraph and $L(X)$ is normal, then
\begin{enumerate}[(a)]
\item $L(X)$ has eigenvalue 0 with multiplicity 1;
\item $\rea(\lambda) > 0$ for $\lambda \neq 0$ an eigenvalue of $L(X)$; and
\item if $L(X) \Zv = \lambda \Zv$, then $L(X)^T \Zv = \comp{\lambda} \Zv$. 
\end{enumerate} \end{lemma}

\proof First, we will show part (c). Let $\Zv$ be an eigenvector of $L := L(X)$ with eigenvalue $\lambda$. Since $L$ and $L^T$ commute, they may be simultaneously diagonalized, so we many assume $\Zv$ is also an eigenvector of $L^T$ with eigenvalue $\theta$. Then
\[
\Zv^*L\Zv = \lambda \Zv^*\Zv
\]
and
\[
\Zv^*L\Zv = (L^*\Zv)^*\Zv = (\theta \Zv)^* \Zv = \comp{\theta} \Zv^*\Zv.
\]
Then, $\comp{\theta} = \lambda$. 

To show parts (a) and (b), we will first show that 
\[
\rea(\lambda) \geq 0
\]
for all eigenvalues $\lambda$ of $L$. Observe that 0 is an eigenvalue of $L$ with the all ones vector as an eigenvector. We will use the property that $X$ is connected to establish that 0 has multiplicity 1 as an eigenvalue and that the other eigenvalues lie in the open right half plane. 

To see that all eigenvalues of $L$ lie in the closed right half plane, we must look at some other incidence matrices. 

For any digraph $X$, we may define the two incidence matrices, $D_h$ and $D_t$ of $X$, with rows indexed by vertices and columns indexed by edges such that the $(u,e)$ entry is as follows:
\[
(D_h)_{u,e} = \begin{cases} 1 & \text{if } u = h(e); \\
                            0 & \text{ otherwise}  \end{cases}
\]
and
\[
(D_t)_{u,e} = \begin{cases} 1 & \text{if } u = t(e); \\
                            0 & \text{ otherwise.}  \end{cases}
\]
We see that if $X$ is eulerian, then 
\[
\Delta^+ = D_tD_t^T = D_h D_h^T
\] and
\[ A = D_tD_h^T. \]
Let $N = D_t - D_h$. Then
\[ \begin{split}
NN^T &= \left(D_t - D_h\right)\left(D_t - D_h\right)^T \\
&= D_t D_t^T - D_t D_h^T - D_hD_t^T + D_hD_h^T\\
&= D_t D_t^T - A + D_hD_h^T  - A^T \\
&= L + L^T.
\end{split} 
\]
But, $NN^T$ is a square, symmetric matrix with all non-negative eigenvalues since 
\[
\Zv^*NN^T\Zv = || N^T \Zv ||^2 \geq 0
\]
for any vector $\Zv \in \cx^n$, where $n$ is the order of $N$. 

Let $\Zv$ be an eigenvector of $L$ with eigenvalue $\lambda$. Then,
\[
(L+L^T)\Zv = (\lambda + \comp{\lambda})\Zv = 2\rea(\lambda) \Zv
\]
Then $2\rea(\lambda)$ is an eigenvalue of $NN^T$ and so $\rea(\lambda) \geq 0$. 

Let $d$ be the maximum out-degree of $A$. Let
\[
B = dI - L 
\]
where $I$ is the $n\times n$ identity matrix. We see that $B$ is a non-negative matrix. Let $\Zv$ be an eigenvector of $L$ with eigenvalue $\lambda$. Then
\[
B\Zv = d\Zv - L\Zv = (d-\lambda)\Zv 
\]
and $\Zv$ is an eigenvector of $B$ with eigenvalue $d-\lambda$. Conversely, if $\Zx$ is an eigenvector of $B$ with eigenvalue $\theta$, then
\[
L\Zx = d\Zx - B\Zx = (d-\theta)\Zx.
\]
Then $\lambda$ is an eigenvalue of $L$ if and only if $d - \lambda$ is an eigenvalue of $B$ with the same multiplicity. 

We may regard $B$ as the adjacency matrix of some digraph $X'$ with some number of loops at each vertex. Since $X$ is connected, then so is $X'$. Then $B$ is irreducible and we may apply the Perron-Frobenius theorem. We obtain that $B$ has a positive real eigenvalue $\rho$ such that $\rho \geq |\theta|$ for all $\theta$ eigenvalues of $B$ and $\rho$ has algebraic multiplicity 1. We see that $\rho = d - \lambda$ for some eigenvalue $\lambda$ of $L$. Since $\rea(\lambda) \geq 0$ for all eigenvalues of $L$, we have that $d -\lambda$ is maximized by $\lambda = 0$. Then $\rho = d$ is the Perron value of $B$ and thus has multiplicity 1. Then 0 is an eigenvalue of $L$ with multiplicity 1. 

If $L$ has an eigenvalue $\lambda = i\beta$ for some $\beta \neq 0 \in \re$. Then $|d - \lambda|$ is an eigenvalue of $B$ and $|d - \lambda| \geq |d|$, which is a contradiction. Then for $\lambda \neq 0$ an eigenvalue of $L$, we see that 
\[
\rea(\lambda) > 0
\]
as claimed. \qed

Observe that if $X$ is not connected, then we may consider the spectrum for each  connected component of $X$ to obtain the following corollary.

\begin{corollary}
If $X$ is a digraph and $L(X)$ is normal, then $\rea(\lambda) >0$ for all eigenvalues $\lambda$ of $L(X)$ except $\lambda = 0$. 
\end{corollary}

\section{Interlacing with the $A$-Laplacian}\label{sec:interlacing}

In this section, we will consider eulerian digraphs $X$ such that the $A$-Laplacian matrix of $X$ is normal. Let $\alpha(X)$ denote the size of the largest acyclic subgraph of $X$. First we will find a spectral bound of disjoint vertex sets $Y$ and $Z$ with no arcs from $Z$ to $Y$. Then, we will use the bound to find a spectral bound for the maximum acyclic subdigraph of $X$. 

We would like to use interlacing to find bounds of combinatorial properties of digraph $X$ using eigenvalues of some matrix respecting the adjacency of $X$. In particular, we would like to use the same method as the proof of Lemma 6.1 in \cite{H95}, which is stated here as Theorem \ref{thm:haemers}. We achieve this using the $A$-Laplacian matrix and prove a result which generalizes the original lemma of Haemers to a sub-class of digraphs. However, as this matrix is not symmetric like in the case for graphs, we need to restrict to digraphs $X$ whose the $A$-Laplacian matrices are normal and we also need a few technical lemmas in order to prove Theorem \ref{thm:offdiag}, the main result of this chapter.

For $L$, we see that the singular values of $L^TL$ are just $|\lambda|$ for each eigenvalue $\lambda$ of $L$.  

\begin{lemma}\label{lem:singperturb} Let $X$ be an eulerian digraph on $n$ vertices such that the Laplacian $L$ of $X$ is normal. If $\lambda$ is an eigenvalue of $L$, then $|\alpha + \lambda|$ is a singular value of 
\[
\alpha I + L
\]
for each $\alpha \in \cx$ and $I$ is the $n\times n$ identity matrix. \end{lemma} 

\proof We need to consider the following:
\[
\begin{split}
(\alpha I + L)(\alpha I + L)^* &= (\alpha I + L)(\comp{\alpha} I + L^T) \\
&= |\alpha|^2 I + \alpha L^T + \comp{\alpha}L + LL^T.
\end{split}
\]
Let $\Zv$ be an eigenvector of $L$ with eigenvalue $\lambda$. Then 
\[
\begin{split}
(\alpha I + L)(\alpha I + L)^*\Zv &=  |\alpha|^2\Zv + \alpha L^T\Zv + \comp{\alpha}L\Zv + LL^T \Zv \\
&=  |\alpha|^2\Zv + \alpha \comp{\lambda} \Zv + \comp{\alpha}\lambda\Zv + |\lambda|^2 \Zv \\
&=  (\alpha + \lambda)(\comp{\alpha + \lambda}) \Zv \\
&= |\alpha + \lambda |^2 \Zv.
\end{split}
\] 
Then the singular values of $\alpha I + L$ are $ |\alpha + \lambda |$ for each eigenvalue $\lambda$ of $L$. \qed 

\begin{lemma} \label{lem:main} Let $X$ be a digraph such that $L(X)$ is normal. Let \[f(\lambda) = \frac{|\lambda|^2}{2 \rea(\lambda)}\] and let $\theta \neq 0$ be the eigenvalue of $L(X)$ which maximizes $f$ amongst non-zero eigenvalues of $L(X)$. Let 
\[ g(\lambda) = \frac{\rea(\lambda)(f(\theta) - f(\lambda))}{\rea(\theta) - \rea(\lambda)}
\]
and let $\mu$ be the eigenvalue of $L(X)$ which minimizes $g$ such that $g(\mu)\geq 0$, if such an eigenvalue exists in the domain of $g$, and $\mu =0$ otherwise.
Let $\widetilde{L} = \alpha I + L(X)$, where $ \alpha \leq - f(\theta)$.
Then $|\alpha|$ is the largest singular value of $\widetilde{L}$.  Further,
\begin{enumerate}[(a)]
\item if $X$ is connected and $\rea(\lambda) \geq \rea(\theta)$ for all $\lambda \notin \{0, \theta\}$, then $|\alpha + \theta|$ is the second largest singular value of $\widetilde{L}$; and
\item if $X$ is connected and $-f(\theta) - g(\mu) \leq \alpha$, then $|\alpha + \theta|$ is the second largest singular value of $\widetilde{L}$.
\end{enumerate}
\end{lemma}

\proof Note that $f$ is well-defined for non-zero eigenvalues $\lambda$ of $L$, since $\rea(\lambda) > 0$ by Lemma \ref{lem:laplace}.  Observe also that $f$ is positive real-valued and $\alpha \in \re$. The function $g$ is well-defined for $\lambda$ when $\rea(\lambda) \neq \rea(\theta)$. If there exists an eigenvalue $\lambda \notin \{0, \theta\}$, such that $\rea(\lambda) < \rea(\theta)$, then, we can see that $g(\lambda) \geq 0$, and so $\mu$ is non-zero. Also, the range for $\alpha$ in part (b), $[-f(\theta) - g(\mu), - f(\theta)]$, is non-empty. 

Since 0 is an eigenvalue of $L$, Lemma \ref{lem:singperturb} gives that $|\alpha|$ is a singular value of $\widetilde{L}$. The singular values of $\widetilde{L}$ are of form $|\alpha + \lambda | $ where $\lambda$ is an eigenvalue of $L$. Let $\lambda$ be a non-zero eigenvalue of $L$. Consider
\begin{equation}\label{eq:main1}
\begin{split}
|\alpha|^2 - |\alpha + \lambda|^2 &= \alpha^2 - (\alpha + \lambda)(\comp{\alpha + \lambda}) \\
&= \alpha^2 - (\alpha + \lambda)(\alpha + \comp{\lambda}) \\
&= \alpha^2 - (\alpha^2 + \alpha(\lambda + \comp{\lambda}) + |\lambda|^2) \\
&= - \alpha(\lambda + \comp{\lambda}) - |\lambda|^2 \\
&= -2\alpha\rea(\lambda) - |\lambda|^2.
\end{split}
\end{equation}
By definition of $\alpha$, we have that
\begin{equation}\label{eq:main2}
- \alpha \geq f(\theta) \geq f(\lambda) = \frac{|\lambda|^2}{2\rea(\lambda)}
\end{equation} 
for all nonzero $\lambda \in \sigma(L)$. From (\ref{eq:main1}) and (\ref{eq:main2}), we obtain:
\[
|\alpha|^2 - |\alpha + \lambda|^2 =  -2\alpha\rea(\lambda) - |\lambda|^2 \geq 0
\]
and we have shown the first part of the statement. 

For statements (a) and (b), let $\delta$ be as follows:
\[
\delta(\lambda) := |\alpha + \theta|^2 - |\alpha + \lambda|^2.
\]
Since $X$ is connected, $L$ has only one eigenvalue whose real part is equal to 0 by Lemma \ref{lem:laplace}. It suffices to show that $\delta(\lambda) \geq 0$ for all nonzero $\lambda \in \sigma(L)$. We expand $\delta(\lambda)$ as follows:
\[
\begin{split}
\delta(\lambda) &= |\alpha + \theta|^2 - |\alpha + \lambda|^2 \\
&= \alpha^2 + \alpha(\theta + \comp{\theta}) + |\theta|^2 - (\alpha^2 + \alpha(\lambda + \comp{\lambda}) + |\lambda|^2)\\
&= |\theta|^2 - |\lambda|^2 + 2\alpha(\rea(\theta) - \rea(\lambda)) .
\end{split}
\]

If $\rea(\theta) = \rea(\lambda)$, then $\delta(\lambda) = |\theta|^2  - |\lambda|^2 $. In this case,
\[
f(\theta) = \frac{|\theta|^2}{2\rea(\theta)} = \frac{|\theta|^2}{2\rea(\lambda)} \geq f(\lambda) = \frac{|\lambda|^2}{2\rea(\lambda)}
\]
and, since $\rea(\lambda)$ and $\rea(\theta)$ are positive, $ |\theta|^2  \geq |\lambda|^2 $ and $\delta(\lambda)\geq 0$. 

If $\rea(\theta) < \rea(\lambda)$, then we may consider
\[ 
\delta(\lambda) = |\theta|^2 - |\lambda|^2 - 2\alpha(\rea(\lambda) - \rea(\theta)) .
\]
Since $(\rea(\lambda) - \rea(\theta)) > 0$ and $- \alpha \geq f(\theta) \geq 0$, we may simplify as follows:
\[ \begin{split}
\delta(\lambda)  &\geq |\theta|^2 - |\lambda|^2 + 2 f(\theta) (\rea(\lambda) - \rea(\theta)) \\
&= |\theta|^2 - |\lambda|^2 +  \frac{|\theta|^2}{\rea(\theta)} (\rea(\lambda) - \rea(\theta)) \\
&= |\theta|^2 - |\lambda|^2 - |\theta|^2 + \frac{|\theta|^2\rea(\lambda)}{\rea(\theta)} \\
&=  - |\lambda|^2 +  \frac{|\theta|^2\rea(\lambda)}{\rea(\theta)} \\
&= \rea(\lambda) \left( - \frac{ |\lambda|^2}{ \rea(\lambda) } +  \frac{|\theta|^2}{\rea(\theta)}\right) \\
&= 2\rea(\lambda) (f(\theta) - f(\lambda)) \\
&\geq 0.
\end{split}
\]
We have shown part (a) and also part (b) when $\rea(\lambda) \geq \rea(\theta)$. 

For part (b), we need only consider eigenvalues $\lambda$ such that $\rea(\theta) > \rea(\lambda)$. In this case, we will use that \[g(\lambda) \geq g(\mu) \geq 0. \] 
Then
\[
\begin{split}
\delta(\lambda) &= |\theta|^2 - |\lambda|^2 + 2\alpha(\rea(\theta) - \rea(\lambda)) \\
&\geq |\theta|^2 - |\lambda|^2 + (- 2f(\theta) -2g(\mu))(\rea(\theta) - \rea(\lambda)) \\
&= |\theta|^2 - |\lambda|^2 - |\theta|^2 + 2f(\theta)\rea(\lambda)  -2g(\mu)(\rea(\theta) - \rea(\lambda)) 
\end{split}
\]
Since $\rea(\theta) - \rea(\lambda) > 0$ and $-2g(\mu) \geq -2g(\lambda)$, we obtain
\[
\begin{split}
\delta(\lambda) &\geq -|\lambda|^2  + 2f(\theta)\rea(\lambda)  -2g(\lambda)(\rea(\theta) - \rea(\lambda)) \\
&=  -|\lambda|^2  + 2f(\theta)\rea(\lambda)  -2 \rea(\lambda)(f(\theta) - f(\lambda)) \\
&= -|\lambda|^2 + 2\rea(\lambda) f(\lambda) \\
&= 0
\end{split}
\]
and we obtain that $\delta(\lambda) \geq 0$, as required. 
 \qed

It is worth observing that if we take $\alpha = -f(\theta)$, then 
\[
|\alpha + \theta|^2 = \alpha^2 + \alpha 2 \rea(\theta) + |\theta|^2 = |\alpha|^2.
\]
We can now prove the main theorem. 

\begin{theorem}\label{thm:offdiag}
Let $X$ be a connected digraph on $n$ vertices where $L := L(X)$ is normal. 
Let $f$, $g$, $\theta$ and $\mu$ be as defined in Lemma \ref{lem:main}. Also, let $\nu \neq 0$ be the eigenvalue of $L$ which minimizes $f$ amongst non-zero eigenvalues of $L$.
Let $Y$ and $Z$ be disjoint vertex sets in $X$ with no arcs from $Z$ to $Y$. Then,
\[
\frac{|Y| |Z|}{(n-|Y|)(n - |Z|)} \leq \frac{ |\alpha + \theta |^2}{\alpha^2} ,
\]
where $ \alpha = \begin{cases} -f(\theta) - f(\nu), & \text{if } \rea(\lambda) \geq \rea(\theta) \text{ for all }\lambda \notin \{0, \theta\}; \\
-f(\theta) - g(\mu), & \text{otherwise.} \end{cases} $ 
% \begin{itemize}
% \item  $ f(\lambda) = \frac{|\lambda|^2}{2 \rea(\lambda)} $,
% \item $\theta \neq 0$ is the eigenvalue of $L$ which maximizes $f$ amongst non-zero eigenvalues of $L$ 
% \item $\nu \neq 0$ is the eigenvalue of $L$ which minimizes $f$ amongst non-zero eigenvalues of $L$ 
% \item $ g(\lambda) = \frac{\rea(\lambda)(f(\theta) - f(\lambda))}{\rea(\theta) - \rea(\lambda)} $, and
% \item $\mu$ is the eigenvalue of $L(X)$ which minimizes $g$ such that $g(\mu)>0$, if such an eigenvalue exists, and $\mu =0$ otherwise.
% \end{itemize}
\end{theorem} 

\proof Let $\alpha =  - f(\theta) - g(\mu)$ and let $\widetilde{L} = \alpha I + L$. In $L$ and $\widetilde{L}$, there is an off-diagonal block of $0$s, where the rows are indexed by $Z$ and columns are indexed by $Y$. This follows directly from hypothesis that there are no arcs from $Z$ to $Y$. We wish to use interlacing to bound the size of such an off-diagonal block of $0$s. Let 
\[
C = \pmat{ 0 & \alpha I + L \\ \alpha I + L^T & 0 } = \pmat{ 0 & \widetilde{L} \\ \widetilde{L}^T & 0 }.
\]
Note that we use $0$ in matrices to represent the zero matrix of the appropriate dimensions and hence will omit subscripts.  
We see that $C$ is symmetric and the eigenvalues of $C$ are 
\[
\{\pm \sigma \,:\, \sigma \text{ a singular value of } \widetilde{L} \}.
\]
Using Lemma \ref{lem:singperturb}, we can write the eigenvalues of $C$ as 
\[
\{\pm |\alpha + \lambda| \,:\, \lambda \text{ eigenvalue value of }  L \}.
\]
By Lemma \ref{lem:main}, we see that $|\alpha|$ is the biggest eigenvalue of $C$ and $|\alpha + \theta|$ is the second largest eigenvalue of $C$.

Since $X$ is eulerian, each row and column of $L$ sums to $0$ and so each row and column of $\widetilde{L}$ sum to $\alpha$. We may partition the rows of $\widetilde{L}$ into rows indexed by $\{ Z,  V(X)\setminus Z\}$ and the columns of $\widetilde{L}$ into columns indexed by $\{ V(X) \setminus Y, Y\}$. Then,
\[
\widetilde{L} = \pmat{ \widetilde{L}_{11} & 0 \\ \widetilde{L}_{21} & \widetilde{L}_{22} }.
\]
This partition of $\widetilde{L}$ induces a partition of $C$ where all diagonal blocks are square;
\[
A = \pmat{0 & 0 & \widetilde{L}_{11} & 0 \\
0 & 0 & \widetilde{L}_{21} & \widetilde{L}_{22} \\
\widetilde{L}_{11}^T & \widetilde{L}_{21}^T & 0 & 0 \\
0 & \widetilde{L}_{22}^T & 0 &0 }.
\]
We let $B$ be the quotient matrix of $C$ with respect to this partition. Recall from Theorem \ref{thm:interlacingquotients} that the entries of $B$ are the average row sums of the corresponding blocks of $C$. We will index the rows and columns of $B$ with $[4]$, for convenience. Since the row and column sums of $\widetilde{L}$ are all equal to $\alpha$, we see that each row and column sum of $\widetilde{L}_{11}$  and of the matrix $\pmat{\widetilde{L}_{21} & \widetilde{L}_{22}}$ is equal to $\alpha$. Then $B(1,3) = B(4,2) = \alpha$ and 
\[
B(2,3) + B(2,4) =  B(3,1) + B(3,2) = \alpha.
\]

Observe that $\widetilde{L}_{22}$ is a $n-z \times y$ matrix, where the lower $y \times y$ submatrix is a principal submatrix of $\widetilde{L}$. For $S,T \subseteq V(X)$, let $E(S,T)$ denote the set of edges $e$ such that $t(e) \in S$ and $h(e) \in T$. Let $W = V(X)\setminus (Y\cup Z)$. We will find $B(2,4)$ by taking the sum over all of the entries of $\widetilde{L}_{22}$ as follows:
\begin{equation}\label{eq:main3}
\begin{split}
(n-2) B(2,4) &= \sum_{j=1}^n \sum_{\ell = 1}^n \widetilde{L}_{22}(j, \ell) \\
&= y\alpha + \sum_{y\ in Y} d^{+}(y) - |E(Y,Y)| - |E(W,Y)|. 
\end{split}
\end{equation} 
Since there are not arcs from $Z$ to $Y$, we have that 
\begin{equation} \label{eq:main4} |E(Y,Y)| + |E(W,Y)| = |E(V(X), Y)| = \sum_{y\in Y} d^{-}(y).
\end{equation} 
Since $X$ is eulerian, we see that $\sum_{y\in Y} d^{-}(y) = \sum_{y\in Y} d^{+}(y)$. Then, from (\ref{eq:main3}) and (\ref{eq:main4}), we obtain that $B(2,4) = \frac{y\alpha}{n-z}$, which implies that $B(2,3) = \alpha - \frac{y\alpha}{n-z}$. 
By an analoguous argument, we find that $B(3,1) = \frac{z\alpha}{n-y}$ and $B(3,2) = \alpha - \frac{z\alpha}{n-y}$. 
Thus we have
\[
B = \pmat{0&0 & \alpha & 0 \\
0 & 0 & \alpha - \frac{y\alpha}{n- z} & \frac{y\alpha}{n- z} \\
\frac{z\alpha}{n- y} & \alpha - \frac{z\alpha}{n- y} & 0 & 0 \\
0 & \alpha & 0 & 0 }.
\]

Since $B$ is the quotient of a symmetric matrix, we see from the proof of Theorem \ref{thm:interlacingquotients} that $B$ is similar to a symmetric matrix. Thus we may let the real eigenvalues of $B$ be 
$\mu_1 \geq \mu_2 \geq \mu _3 \geq \mu_4 $.
Let 
$\lambda_1 \geq \lambda_2 \geq \cdots \geq \lambda_{2n-1} \geq \lambda_{2n}$
be the eigenvalues of $C$. Observe that $C$ is similar to $-C$ and $B$ is similar to $-B$ by construction of $C$. Then, we have that 
\[
\mu_4 = - \mu_1, \, \mu_3 = -\mu_2, \, \lambda_{2n} = - \lambda_1 \text{ and } \lambda_{2n-1} = - \lambda_2. 
\]
Applying the interlacing theorem gives
\[
\lambda_1 \geq \mu_1 ,\, \lambda_2 \geq \mu_2
\text{ and }
\mu_3 \geq \lambda_{2n-1} ,\, \mu_4 \geq \lambda_{2n}.
\]
Recall that $(\lambda_1, \lambda_2) = ( |\alpha| , |\alpha + \theta |) $. Then 
\begin{equation}\label{eq:main5}
\mu_1\mu_2 \mu_3\mu_4 = (-1)^2(\mu_1\mu_2)^2 \leq (\lambda_1\lambda_2)^2 =\left( |\alpha| |\alpha + \theta | \right)^2.
\end{equation}
On the other hand, we see that 
\begin{equation}\label{eq:main6}
\mu_1\mu_2 \mu_3\mu_4 = \det(B) =  \frac{\alpha^2 y}{n-z} \frac{\alpha^2 z}{n-y}.
\end{equation}
From (\ref{eq:main5}) and (\ref{eq:main6}), we obtain that 
\[
  \frac{\alpha^2 y}{n-z} \frac{\alpha^2 z}{n-y} \leq ( |\alpha| |\alpha + \theta |)^2
\]
which simplifies to 
\[
 \frac{y}{n-z} \frac{z}{n-y} \leq \frac{ |\alpha + \theta |^2}{\alpha^2}
\]
as claimed. \qed

Observe that if $X$ is a graph, then the $A$-Laplacian is the usual Laplacian matrix of a graph. In this case, the $A$-Laplacian is symmetric and hence normal and so all eigenvalues are real and non-negative. Thus, for $\lambda$ an eigenvalue of $L(X)$, we see that 
\[
f(\lambda) = g(\lambda) = \lambda/2.
\]
Then, we can recover the original theorem of Haemers \cite[Lemma 6.1]{H95} for graphs as a corollary of Theorem \ref{thm:offdiag}.

A tournament has normal adjacency matrix if and only if it is regular (see \cite{dC91}). Then, the $A$-Laplacian matrices of regular tournaments are normal matrices. Let $X$ be a regular tournament on $n$ vertices. In this case, all of the non-zero eigenvalues have real part equal to $\frac{n}{2}$. We see that 
\[ f(\lambda) = \frac{|\lambda|^2}{2\frac{n}{2}} =\frac{|\lambda|^2}{n}. \] 
Then, $\alpha = - \frac{|\theta|^2}{n} -  \frac{|\nu|^2}{n}$, where $\theta$ and $\nu$ are the largest and smallest eigenvalues of $L(X)$ in magnitude. Theorem \ref{thm:offdiag} gives the following corollary.

\begin{corollary}\label{cor:tournaments}
 Let $X$ be a regular tournament on $n$ vertices and $(Z, Y)$ be a separation in $X$. Then
\[
\frac{|Y| |Z|}{(n-|Y|)(n - |Z|)} \leq \frac{ |\alpha + \theta |^2}{\alpha^2} ,
\]
where $ \alpha = - \frac{|\theta|^2}{n} -  \frac{|\nu|^2}{n}$ where $\theta$ and $\nu$ are the largest and smallest non-zero eigenvalues of $L(X)$ in magnitude.
\end{corollary}

%\bibliographystyle{plain}
%\bibliography{spectra}

\begin{thebibliography}{10}

\bibitem{BePl94}
Abraham Berman and Robert~J. Plemmons.
\newblock {\em Nonnegative matrices in the mathematical sciences}, volume~9 of
  {\em Classics in Applied Mathematics}.
\newblock Society for Industrial and Applied Mathematics (SIAM), Philadelphia,
  PA, 1994.
\newblock Revised reprint of the 1979 original.

\bibitem{CG97}
Ada Chan and Chris~D. Godsil.
\newblock Symmetry and eigenvectors.
\newblock In {\em Graph symmetry ({M}ontreal, {PQ}, 1996)}, volume 497 of {\em
  NATO Adv. Sci. Inst. Ser. C Math. Phys. Sci.}, pages 75--106. Kluwer Acad.
  Publ., Dordrecht, 1997.

\bibitem{C89}
F.~R.~K. Chung.
\newblock Diameters and eigenvalues.
\newblock {\em J. Amer. Math. Soc.}, 2(2):187--196, 1989.

\bibitem{Cv72}
Drago{\v{s}}~M. Cvetkovi{\'c}.
\newblock Chromatic number and the spectrum of a graph.
\newblock {\em Publ. Inst. Math. (Beograd) (N.S.)}, 14(28):25--38, 1972.

\bibitem{dC91}
D.~de~Caen.
\newblock The ranks of tournament matrices.
\newblock {\em Amer. Math. Monthly}, 98(9):829--831, 1991.

\bibitem{H95}
Willem~H. Haemers.
\newblock Interlacing eigenvalues and graphs.
\newblock {\em Linear Algebra Appl.}, 226/228:593--616, 1995.

\bibitem{H70}
Alan~J. Hoffman.
\newblock On eigenvalues and colorings of graphs.
\newblock In {\em Graph {T}heory and its {A}pplications ({P}roc. {A}dvanced
  {S}em., {M}ath. {R}esearch {C}enter, {U}niv. of {W}isconsin, {M}adison,
  {W}is., 1969)}, pages 79--91. Academic Press, New York, 1970.

\bibitem{M91}
Bojan Mohar.
\newblock Eigenvalues, diameter, and mean distance in graphs.
\newblock {\em Graphs Combin.}, 7(1):53--64, 1991.

\bibitem{M10}
Bojan Mohar.
\newblock Eigenvalues and colorings of digraphs.
\newblock {\em Linear Algebra Appl.}, 432(9):2273--2277, 2010.

\bibitem{MoPo93}
Bojan Mohar and Svatopluk Poljak.
\newblock Eigenvalues in combinatorial optimization.
\newblock In {\em Combinatorial and graph-theoretical problems in linear
  algebra ({M}inneapolis, {MN}, 1991)}, volume~50 of {\em IMA Vol. Math.
  Appl.}, pages 107--151. Springer, New York, 1993.

\bibitem{R07}
Peter Rowlinson.
\newblock The main eigenvalues of a graph: a survey.
\newblock {\em Appl. Anal. Discrete Math.}, 1(2):445--471, 2007.

\bibitem{R13}
Peter Rowlinson.
\newblock On graphs with an eigenvalue of maximal multiplicity.
\newblock {\em Discrete Math.}, 313(11):1162--1166, 2013.

\end{thebibliography}

\end{document}